\documentstyle[mssymb]{article}
\newcommand{\be}{\begin{equation}}
\newcommand{\ee}{\end{equation}}
\newcommand{\ra}{\rightarrow}
\newcommand{\ul}{\underline}

\newcommand{\R}{{\bf R}}
\newcommand{\Z}{{\bf Z}}
\font\svnbf=cmbx7
\newcommand{\sZ}{{\hbox{\svnbf Z}}} % small \Z for subscripts
\font\bigbf=cmbx12 scaled \magstep1
\newcommand{\lZ}{{\hbox{\bigbf Z}}} % large \Z for title
\begin{document}
\centerline{\Large A characterization of the $\lZ^n$ lattice}
\vspace*{3ex}
\centerline{Noam D. Elkies}
\vspace*{6ex}

\noindent{\bf Introduction.}
\vspace*{1ex}

In this note we prove that $\Z^n$ is the only integral unimodular
lattice $L\subset\R^n$ which does not contain a vector $w$ such that
$|w|^2<n$ and $(v,v+w)\equiv0\bmod2$ for all $v\in L$.  By the work
of Kronheimer and others on the Seiberg-Witten equation this yields
an alternative proof of a theorem of Donaldson~[D1,D2]
on the geometry of {4-manifolds}.

The proof uses the theory of theta series and modular forms;
since this technique is not yet in the standard-issue arsenal
of the {4-manifold} community, I begin with an abbreviated
exposition of this theory to make this note reasonably
self-contained.  This develops only the barest minimum,
even to the point of never using the phrase ``modular form'';
for a more substantial exposition, refer to \cite[Ch.VII]{Serre},
and note the concluding remarks (6.7, ``Complements'').

Knowing that any $L\not\cong\Z^n$ has characteristic vectors
of norm $\leqslant n-8$, one might ask for which lattices is
$n-8$ the minimum.  It turns out that the same technique
also yields a complete answer to this question.  Since the
answer may be of some interest (for instance there are 14
such lattices in each dimension $n\geqslant23$), but its proof
requires a somewhat more extensive use of modular forms,
we announce the result at the end of this note but defer
its proof and further discussion to a later paper.

\vspace*{3ex}

\noindent{\bf Fractional linear transformations and theta series.}

\vspace*{1ex}

Let $H$\/ be the Poincar\'e upper half-plane $\{ t = x+iy : y>0 \}$,
and let $\Gamma$ be the group PSL$_2(\Z)={\rm SL}_2(\Z)/\{\pm{\bf1}\}$,
acting on~$H$\/ by fractional linear transformations:
\be
\left({a\;b\atop c\;d}\right) : \ t \mapsto \frac{at+b}{ct+d} \; .
\ee
It is known that $\Gamma$ is generated by
$S=({0\;-1\atop1\;\phantom-0})$ and $T=({1\;1\atop0\;1})$,
acting on~$H$\/ by
\be
S(t) = - \frac 1 t, \quad T(t) = t + 1.
\ee
Let $\Gamma(2)=\{({a\;b\atop c\;d})\in\Gamma:b,c{\rm\ even}\}$;
this is a normal subgroup of~$\Gamma$,
and reduction mod~2 yields the quotient map $\Gamma\ra\Gamma/\Gamma(2)
={\rm PSL}_2(\Z/2)\cong{\rm S}_3$.  Finally let $\Gamma_+\subset\Gamma$
be the subgroup generated by~$S$\/ and~$T^2$.  Then $\Gamma_+$ has
index~3 in~$\Gamma$, contains $\Gamma(2)$ with index~2, and consists
of the matrices congruent mod~2 to either {\bf1} or~$S$.  Indeed it
is clear that $\Gamma_+$ consists of matrices of this form;
that all such matrices are in~$\Gamma_+$ is perhaps most readily seen
by proving as in \cite[Ch.VII,~Thm.1,2]{Serre} that
\be
D_+ := \{ t = x+iy \in H : |x|\leq1, |t|\geq1 \}
\ee
(the ideal hyperbolic triangle in~$H$\/ with vertices $-1,1,i\infty$)
is a fundamental domain for the action of $\Gamma_+$ on~$H$,
and noting that $D_+$ is 3 times as large as the standard
fundamental domain for~$\Gamma$.

Now let $L$ be a unimodular integral lattice in $\R^n$,
i.e.\ a lattice of discriminant~1 such that $(v,v')\in\Z$
for all $v,v'\in L$.  The \ul{theta} \ul{series}
$\theta_L$ of~$L$ is a generating function encoding the norms
$|v|^2=(v,v)$ of lattice vectors:
\be
\theta_L(t) := \sum_{v\in L} e^{\pi i |v|^2 t} \quad (t\in H).
\label{thetadef}
\ee
For instance for $n=1$ we have
\be
\theta_\sZ(t) := 1 +
2 \left( e^{\pi i t} + e^{4\pi i t} + e^{9\pi i t} + \cdots \right) \; .
\label{thetaZ}
\ee
This sum converges uniformly in compact subsets of~$H$\/
(if $t=x+iy$ then $|e^{\pi i |v|^2 t}| = e^{-\pi|v|^2 y}$)
and thus defines a holomorphic function on~$H$.  If $L_1,L_2$
are unimodular integral lattices in $\R^{n_1},\R^{n_2}$ then
$L_1\oplus L_2$ is a unimodular integral lattice in $\R^{n_1+n_2}$
whose theta series is given by
\be
\theta_{L_1\oplus L_2}(t) =
\theta_{L_1}(t) \cdot \theta_{L_2}(t).
\label{thetaprod}
\ee
Since each $|v|^2$ is an integer we have
\be
\theta_L(t) = \theta_L(t+2) = \theta_L(T^2(t)).
\label{thetaT}
\ee
Since $L$ is its own dual lattice we obtain a more interesting
functional equation by applying Poisson inversion to~(\ref{thetadef}):
\be
(t/i)^{n/2} \theta_L(t) = \theta_L(-1/t) = \theta_L(S(t)),
\label{thetaS}
\ee
where $(t/i)^{n/2}$ is the $n$th power of the principal branch
of $\sqrt{t/i}$.  By iterating (\ref{thetaT},\ref{thetaS}) we
find that for every $g=({a\;b\atop c\;d})$ in
$\langle S,T^2 \rangle = \Gamma_+$ there is a functional equation
\be
\theta_L(g(t)) = \epsilon_n(c,d) \cdot (ct+d)^{n/2} \theta_L(t),
\label{thetaG}
\ee
where again $(ct+d)^{n/2}$ is the $n$th power of the principal branch
of $\sqrt{ct+d}$, and $\epsilon_n(c,d)$ is an eighth root of unity
which does not depend on the choice of unimodular integral lattice~$L$.
(It does not depend on $a,b$ because $c,d$\/ determine $g$ up to
a power of~$T^2$.)  By choosing $L=\Z^n$ and using (\ref{thetaprod})
we find
\be
\epsilon_n(c,d) = \bigl( \epsilon_1(c,d) \bigr)^n.
\label{epsilon}
\ee
Note that \cite[Ch.VII]{Serre} assumes that $L$\/ is an \ul{even}
lattice, i.e.\ $|v|^2\in 2\Z$ for all $v\in L$.  Such~$L$ have theta
series invariant under~$T$, and thus satisfy (\ref{thetaG}) for all
$g\in\langle S,T\rangle=\Gamma$.  It is known from the arithmetic
theory \cite[Ch.V]{Serre} that $n\equiv0\bmod8$ for such lattices, 
whence the $\epsilon_n$ factors all equal~1 in that case; this
could also be proved analytically using (\ref{thetaS}) and the
identity $(ST)^3=\bf1$.  We shall soon observe, en route to our
estimate on the norm of characteristic vectors of odd lattices,
that this method also yields an analytic proof of the fact
\cite[Ch.V,~Thm.2]{Serre} that these vectors all have norm
$\equiv n \bmod 8$.

How do fractional linear transformations $g\in\Gamma-\Gamma_+$
act on $\theta_L$?  We need only consider one representative
of each of the two nontrivial cosets of $\Gamma_+$ in $\Gamma$,
for instance $g=T$ and $g=TS$.  For the first we find simply
\be
\theta_L(T(t)) = \theta_L(t+1) = \sum_{v\in L} e^{\pi i |v|^2 (t+1)}
= \sum_{v\in L} (-1)^{|v|^2} e^{\pi i |v|^2 t} .
\label{thetat+1}
\ee
Now recall that the sign $v\mapsto(-1)^{|v|^2}$ is a group homomorphism
from~$L$ to $\{\pm1\}$ (because
\be
|v+v'|^2 = |v|^2 + |v'|^2 + 2(v,v') \equiv |v|^2 + |v'|^2 \bmod 2
\label{mod2}
\ee
for all $v,v'\in L$).  Since $L$ is unimodular there is a bijection
between characters $L\ra\{\pm1\}$ and cosets of~$2L$ in~$L$
which associates to the coset of any $w\in L$ the character
$v\mapsto(-1)^{(v,w)}$.  In particular there is a coset associated
with $v\mapsto(-1)^{|v|^2}$; vectors in that coset, characterized by
\be
|v|^2 \equiv (v,w) \bmod 2 {\rm\ for\ all}\ v\in L,
\label{chardef}
\ee
are known as \ul{characteristic} \ul{vectors} of~$L$.
(In \cite[Ch.V]{Serre} this coset is called the ``canonical class''
in $L/2L$; in~\cite{IEEE} this coset, scaled by~$1/2$ to obtain a
translate of~$L$ by $w/2$, is called the ``shadow'' of~$L$, and
our key formula (\ref{thetaTS}) below is also a key ingredient
of~\cite{IEEE}.)  Choose some characteristic vector~$w$, and rewrite
(\ref{thetat+1}) as
\be
\theta_L(t+1) =
\sum_{v\in L} e^{\pi i \left(|v|^2 t + (v,w)\right)} .
\label{thetat+1'}
\ee
Applying Poisson inversion to this sum we find
\be
(t/i)^{n/2}\theta_L(t+1) =
\sum_{v\in L} e^{\pi i |v+\frac{w}{2}|^2(\frac{-1}{t})}
= \theta'_L(S(t)),
\label{Pois2}
\ee
where 
\be
\theta'_L(t) := \sum_{v\in L+\frac{w}{2}} e^{\pi i |v|^2 t}
\label{theta'def}
\ee
is a generating function encoding the norms of canonical vectors.
Replacing $t$ by $St=-1/t$ in~(\ref{theta'def}) we conclude that
\be
\theta_L(TS(t)) = \theta_L(\frac{-1}{t}+1) = (t/i)^{n/2}\theta'_L(t).
\label{thetaTS}
\ee
To recover the result
\be
|w|^2\equiv n\bmod 8
\label{mod8}
\ee
we may now regard (\ref{thetaTS}) as a formula for $\theta'_L(t)$
and compare it with
\samepage{
$$
\left(\frac{t+1}{i}\right)^{\!n/2}\theta'_L(t+1) =
\theta_L(TST(t)) = \theta_L(ST^{-1}S(t))
$$ \vspace*{-5ex} \be \phantom. \label{ST3} \ee \vspace*{-5ex} $$
= (T^{-1}S(t)/i)^{n/2} \theta_L(T^{-1}S(t)) =
\left(\frac{i(t+1)}{t}\right)^{\!n/2} \theta_L(TS(t))
%\label{ST3}
$$
}
\vspace*{-2ex} \\ %grrr
(in which we used $S^2=(ST)^3=\bf1$ and the invariance of~$\theta_L$
under $T^2$, and again use $n/2$ power to mean $n$th power of principal
square root).  This yields
\be
\theta'_L(t+1) = e^{\pi i n/4} \theta'_L(t).
\label{aha8}
\ee
Thus $\theta'_L(t)$ is a linear combination of terms $e^{\pi i m t/4}$
with $m\equiv n\bmod8$, from which it follows that all the
characteristic vectors have norm $\equiv n\bmod8$ as claimed.

The characteristic vectors of~$\Z$ are the odd integers, so
\be
\theta'_\sZ(t) = 2\sum_{m=0}^\infty e^{\pi i (m+\frac12)^2 t} =
2 e^{\pi i t/4} \left(
1 + e^{2\pi i t} + e^{6\pi i t} + e^{12\pi i t} + \cdots \right) \; .
\label{theta'Z}
\ee
Thus $\theta'_\sZ(t)\sim2e^{\pi i t/4}\ra0$ as $t\ra i\infty$.
From~(\ref{thetaTS}) it follows that $\theta_\sZ$ tends to~zero
as $t\in D_+$ approaches the ``cusp'' $\pm1$.  It will be crucial
to us that $\theta_\sZ$ {\bf has no zeros in~$H$}.  This can be seen
either from explicit product formulas such as
\be
\sum_{m=0}^\infty q^{(m+\frac12)^2} =
q^{1/4} \prod_{j=1}^\infty (1+q^{2j})(1-q^{4j})
\label{jacobiprod}
\ee
(a special case of the Jacobi triple product),
or by using contour integrals as in~\cite[Ch.VII,~Thm.3]{Serre}
to show that $\pm1$ is the only zero of~$\theta_\sZ$ in
$D_+\cup\{{\rm cusps}\}$.  Also $\theta_\sZ(i\infty)=1$ so
$\theta_\sZ$ is bounded away from zero as $t\ra i\infty$.

\vspace*{3ex}

\noindent{\bf The shortest characteristic vector.}

\vspace*{1ex}

We are now ready to prove:

{\bf Theorem.} {\sl Let $L$ be a unimodular integral lattice
in $\R^n$ with no characteristic vector $w$ such that $|w|^2<n$.
Then $L\cong\Z^n$.
}

{\sl Proof}\/: We first show that $L$ and $\Z^n$ have the same
theta function.  To that end consider
\be
R(t) := \theta_L(t) / \theta_{\sZ^n}(t)
= \theta_L(t) / \theta^n_\sZ(t).
\label{Rdef}
\ee
This is a holomorphic function because $\theta_\sZ$ does not vanish
in~$H$.  Since $\theta_L$ and $\theta_{\sZ^n}$ both transform according
to~(\ref{thetaG}) under~$\Gamma_+$, their quotient $R(t)$ is invariant
under~$\Gamma_+$.  By the hypothesis on~$L$ we have
$\theta'_L\ll e^{\pi i n t/4}$ as $t\ra i\infty$.  Thus
$\theta'_L/\theta'_{\sZ^n}$ is bounded as $t\ra i\infty$,
whence by~(\ref{thetaTS}) $R(t)$ is bounded as $t\in D_+$
approaches $\pm 1$.  Finally $R(i\infty)=1$.  By the maximum
principle we deduce that $R$\/ is the constant function~1,
i.e.\ $\theta_L=\theta_\sZ^n$.

Thus for each $m$ the lattices $L$ and $\Z^n$ have the same number
of vectors of norm~$m$.  Taking $m=1$ we find that $L$ has
$n$ pairs of unit vectors.  Since $L$ is integral these
must be orthogonal to each other, and thus generate a copy
of~$\Z^n$ inside~$L$.  Using integrality again we conclude that
this copy is all of~$L$.~~QED

\vspace*{1ex}

Since the hypothesis is automatically satisfied if $n<8$ we also
recover the fact that $\Z^n$ is the only unimodular integral lattice
for those~$n$.  With some more work we can also use the relation
between $\theta_L$ and $\theta'_L$ and the theory of modular forms to
completely describe those $L\subset\R^n$ whose shortest characteristic
vector has norm $n-8$: these are precisely the lattices of the form
\hbox{$\Z^{n-r} \oplus L_0$}, where $L_0\subset\R^r$ is a unimodular
integral lattice with no vectors of norm~1 and exactly \hbox{$2n(23-n)$}
vectors of norm~2.  In particular $n\leqslant23$ and there are
only finitely many choices for $L_0$.  Fortunately the table
of unimodular lattices in~\cite[pp.416--7]{SPLAG} extends just
far enough that we can list all possible $L_0$.  These are tabulated
below, indexed as in the table of~\cite{SPLAG} by the root system
of {norm-2} vectors:

\vspace*{1ex}

\centerline{
$
\begin{array}{c|cccccccccccccc}
r & 8 & 12 & 14 & 15 & 16 & 17 & 18 & 18 & 19 & 20 & 20 & 21 & 22 & 23
\\ \hline
\phantom{\sum^0}\!\!\!\!\! & E_8 & D_{12} & E_7^2 & A_{15} & D_8^2 & A_{11}E_6
 & D_6^3 & A_9^2 & A_7^2D_5 & D_4^5 & A_5^4 & A_3^7 & A_1^{22} & O_{23}
\end{array}
$
}

\vspace*{1ex}

Of these the first is the $E_8$ lattice and the last is the
``shorter Leech lattice''; these are the unimodular integral
lattices of minimal dimension having minimal norm~2 and~3
respectively.  It also follows from the analysis that each
of these lattices has exactly $2^{n-11} r$ characteristic vectors
of norm~$n-8$.  We defer the proof of the $\Z^{n-r}\oplus L_0$
criterion, and an analogous condition for self-dual binary codes,
to a subsequent paper.

{\bf Acknowledgements.}  Thanks to Tom Mrowka for bringing this
problem to my attention and to John H. Conway for enlightening
correspondence.  This work was made possible in part by funding from
the National Science Foundation and the Packard Foundation.

%\nopagebreak

\vspace*{5ex}
\begin{small}
\noindent 
Dept.\ of Mathematics\\
Harvard University\\
Cambridge, MA 02138 USA\\
e-mail: {\tt elkies@math.harvard.edu}\\ \\
January, 1995
\end{small}

\end{document}